\documentclass{amsart}
\usepackage{graphicx}
\newtheorem{thm}{Theorem}[section]
\newtheorem{cor}[thm]{Corollary}
\newtheorem{lem}[thm]{Lemma}

\theoremstyle{definition}

\theoremstyle{remark}
\newtheorem{rem}[thm]{Remark}
\numberwithin{equation}{section}
\begin{document}

\title[]{Multiplicative derivations on rank-$S$ matrices for relatively small $S$}%
\author{Xiaowei Xu}%
\address{College of Mathematics, Jilin University, Changchun 130012, China}%
\email{xuxw@jlu.edu.cn}%

\author{Baochuan Xie}%
\address{College of Mathematics, Jilin University, Changchun 130012, China}%
\email{2956015407@qq.com}%

\author{Yanhua Wang}%
\address{School of Mathematics, Shanghai Key Laboratory of Financial Information Technology, Shanghai University of Finance and Economics, Shanghai 200433,China}%
\email{yhw@mail.shufe.edu.cn}%

\author{Zhibing Zhao}%
\address{School of Mathematical Sciences, Anhui University, Hefei 230601, China}%
\email{zbzhao@ahu.edu.cn}%

\thanks{}%

%\date{}%
%\dedicatory{}%
%\commby{}%
% ----------------------------------------------------------------
\begin{abstract}
Let $n$ and $s$ be fixed integers such that $n\geq 2$ and $1\leq s\leq \frac{n}{2}$. Let $M_n(\mathbb{K})$ be the ring of all $n\times n$ matrices over a field $\mathbb{K}$. If  a map
$\delta:M_n(\mathbb{K})\rightarrow M_n(\mathbb{K})$ satisfies that $\delta(xy)=\delta(x)y+x\delta(y)$ for any two rank-$s$ matrices $x,y\in M_n(\mathbb{K})$, then there exists a derivation $D$ of
$M_n(\mathbb{K})$ such that $\delta(x)=D(x)$  holds for each rank-$k$ matrix $x\in M_n(\mathbb{K})$ with $0\leq k\leq s$.
\end{abstract}
\maketitle
% ----------------------------------------------------------------
\section{Introduction}
Franca \cite{Franca-2012-LAA} initialed the research on nonadditive subsets of prime rings in the theory of functional identities   by describing the commuting additive map
on the set of all $n\times n$ invertible matrices or the set of  all $n\times n$ singular matrices  rather than the ring of all $n\times n$
matrices over fields.
This is an extension of the well-known theorem
of Bre\v{s}ar (see the original paper \cite[Theorem A]{Bresar-1993-JAlg}, or the survey paper \cite[Corollary
3.3]{Bresar-2004-Taiwan},  or the book \cite[Corollary5.28]{Bresar-Chebotar-Martindale-2007-JAlg}).
Furthermore, in 2013, Franca \cite{Franca-2013-LAA} (also see Xu et al. \cite{XuYi-2014-ToELA})
 extended the discussion  to the set of all rank-$s$ matrices over fields for fixed $2\leq s<n$.
In 2014, Liu (see
\cite{LiuCK-centralizing,LiuCK-strong}) researched centralizing additive maps and strong commutativity
preserving maps on the set of all $n\times n$ invertible matrices or  the set of all $n\times n$ singular matrices over division
rings and obtained nice
conclusions, which developed the corresponding results in the theory of functional identities.
Recently, Xu et al. \cite{XuLiu-LAMA,XuYi-2014-ToELA} proved that
a map $g$ from the ring of all $n\times n$
matrices over a field into itself is additive if and only if $g(A+B)=g(A)+g(B)$ for any two rank-$s$ matrices $A,B
\in M_n(\mathbb{K})$, where $\frac{n}{2}\leq s\leq n$ is fixed. For  further references see \cite{XuLiZhu-2017-Aequat-Math,Kosan-Sahinkaya-Zhou-2018-CMB,Liu-Liau-Tsai-2018-OAM,Liu-Yang-2017-LAA,Franca-2017-LAMA,XuZhu-2018-LAMA}.

 On the other hand,  a map $f$ from a ring $R$ into itself is called a multiplicative isomorphism if
 $f$ is  bijective and $f(xy)=f(x)f(y)$ for all $x,y\in R$.  A map $f$ from a ring $R$ into itself is called a multiplicative  derivation if
 $f(xy)=f(x)y+xf(y)$ for all $x,y\in R$.
The question of when a multiplicative
isomorphism is additive has been considered by Rickart \cite{Rickart-1948-BAMS} and Johnson \cite{Johnson-1958-PAMS}.
Martindale \cite{Martindale-1969-PAMS} improved the main theorem of Rickart \cite[Theorem II]{Rickart-1948-BAMS}.

In 1991, Daif \cite{Daif-1991-IJMS} considered the similar question of when a multiplicative derivation is
additive. He proved that it is true for the ring $R$ with an idempotent element $e\neq 0,1$ satisfying: (1) $xR=0$ implies $x=0$;
(2) $eRx=0$ implies $x=0$;
(3) $exeR(1-e)=0$ implies $exe=0$. Note that for $n\geq 2$ the ring $M_n(R)$ ($T_n(R)$, respectively) of all $n\times n$  (upper triangular) matrices over a unital ring $R$ is a special example of the rings Daif stated. So a multiplicative derivation of $M_n(R)$ ($T_n(R)$, respectively) must be a derivation,
where $R$ is a ring with an identity and $n>1$.

In this short note, we consider the multiplicative derivation on the set of all $n\times n$ rank-$s$ matrices over a field $\mathbb{K}$ other than
 the ring of all  $n\times n$ matrices over $\mathbb{K}$ and prove that for the case $1\leq s\leq \frac{n}{2}$  and $n\geq 2$, if a map $\delta:M_n(\mathbb{K})\rightarrow M_n(\mathbb{K})$ satisfies that $\delta(xy)=\delta(x)y+x\delta(y)$ for any two rank-$s$ matrices $x,y\in M_n(\mathbb{K})$, then there exists a derivation $D$ of
$M_n(\mathbb{K})$ such that $\delta(x)=D(x)$  for each rank-$k$ matrix $x\in M_n(\mathbb{K})$ with $k\leq s$.
This means that the multiplicative derivation on rank-$s$ matrices over a field is almost a derivation when restricted on the matrices whose rank  is not more than $s$ for relative small $s$. As an application, we will show that the multiplicative derivation on some nonadditive subset of the matrix ring $M_n(\mathbb{K})$ over a field $\mathbb{K}$ has to be a derivation.

\section{Multiplicative derivations on rank-$s$ matrices for relatively small $s$}

In this section, unless stated otherwise, we will always assume that both  $n$ and $s$ are fixed integers such that $n\geq 2$ and $1\leq s\leq\frac{n}{2}$, and always  denote by $\mathbb{K}$ a
 field, by $M_n(\mathbb{K})$ the ring of all $n\times n$
matrices over $\mathbb{K}$, by $GL_n(\mathbb{K})$  the set of all $n\times n$ invertible matrices over $\mathbb{K}$.
For $0\leq k\leq n$, the symbol  $M_n^{k}(\mathbb{K})$ ($M_n^{\leq k}(\mathbb{K})$ and $M_n^{<k}(\mathbb{K})$, respectively) will always
 denote the set of all matrices whose rank  is equal to (not more than and less than, respectively) $k$ in
$M_n(\mathbb{K})$. A map
$\delta: M_n(\mathbb{K})\rightarrow M_n(\mathbb{K})$ is called a {\it multiplicative derivation on $\mathcal{S}$ a subset of $M_n(\mathbb{K})$ if $\delta(xy)=\delta(x)y+x\delta(y)$ for all $x,y\in \mathcal{S}$. Write $\mathcal{D}_\mathcal{S}^{\times}(M_n(\mathbb{K}))$ for the set of all
multiplicative derivations on the subset $\mathcal{S}$ of $M_n(\mathbb{K})$. If $\mathcal{S}=M_n^s(\mathbb{K})$, we also write $\mathcal{D}_s^{\times}(M_n(\mathbb{K}))$ for $\mathcal{D}_\mathcal{S}^{\times}(M_n(\mathbb{K}))$ and call a multiplicative derivation on $M_n^s(\mathbb{K})$
a multiplicative derivation on rank-$s$ matrices. Write $e_{ij}$ for the $n\times n$ matrix with 1 in the position
$(i,j)$ and 0 in every other position. The symbol $\sum_{i=a}^be_{ii}$ will denote zero matrix once $a>b$. Denote by $I_n$  the $n\times n$ identity matrix, by $\underline{n}$ the set $\{1,2,\ldots,n\}$ and by
$\mathbb{K}^t$ the set of all $t\times 1$ matrices over $\mathbb{K}$.

Firstly, we note that the set of all multiplicative derivations on a nonempty subset $\mathcal{S}$ of  $M_n(\mathbb{K})$ is a vector space.
\begin{lem}
 $\mathcal{D}^{\times}_\mathcal{S}(M_n(\mathbb{K}))$ is a $\mathbb{K}$-vector space. \label{lemma-vector-space}
\end{lem}

{\bf Proof.} We only need to show that for any $\delta_1,\delta_2\in \mathcal{D}_\mathcal{S}^{\times}(M_n(\mathbb{K}))$ and any $\lambda_1,\lambda_2\in \mathbb{K}$
\begin{equation}\label{equation-lemma-vector-space}
\lambda_1\delta_1+\lambda_2\delta_2\in \mathcal{D}_\mathcal{S}^{\times}(M_n(\mathbb{K})).
\end{equation}
In fact, for any $x,y\in \mathcal{S}$,
$$
\begin{array}{cl}
&(\lambda_1\delta_1+\lambda_2\delta_2)(xy)\\
=&(\lambda_1\delta_1)(xy)+(\lambda_2\delta_2)(xy)
=\lambda_1\delta_1(xy)+\lambda_2\delta_2(xy)\\
=&\lambda_1\delta_1(x)y+\lambda_1x\delta_1(y)+\lambda_2\delta_2(x)y+\lambda_2x\delta_2(y)\\
=&(\lambda_1\delta_1+\lambda_2\delta_2)(x)y+x(\lambda_1\delta_1+\lambda_2\delta_2)(y)
\end{array}
$$
which implies that \eqref{equation-lemma-vector-space} holds. \hfill $\Box$

The following Remark \ref{remark-y-y1-y2} and  Corollary \ref{corollary-y-y1-y2} will be used in the proof of Lemma \ref{lemma-main}, Theorem \ref{Theorem-rank-s} and Corollary \ref{corolary-log}.
\begin{rem}
Let $M_n(\mathbb{K})$ be the ring of all $n\times n$
matrices over a field $\mathbb{K}$ where $n\geq 2$. Let $1\leq s\leq n$ and $2s-n\leq k\leq s$ be integers.   Then for each $y\in M_n^{k}({\mathbb{K}})$, there exist
$y_1,y_2\in M_n^{s}({\mathbb{K}})$ such that $y=y_1y_2$. \label{remark-y-y1-y2}
\end{rem}

{\bf Proof.} There exist invertible matrices $P,Q\in GL_n(\mathbb{K})$ such that $$y=P\left(\sum_{i=1}^ke_{ii}\right)Q,$$
where we denote  by $\sum_{i=1}^ke_{ii}$ the zero matrix
in the case of $k=0$. From $s+s-k=2s-k\leq n$ we have the desired matrices
$$y_1=P\left(\sum_{i=1}^ke_{ii}+\sum_{j=k+1}^{s}e_{jj}\right),~~y_2=\left(\sum_{i=1}^ke_{ii}+\sum_{j=s+1}^{2s-k}e_{jj}\right)Q.
$$ \hfill$\Box$

\begin{cor}
Let $M_n(\mathbb{K})$ be the ring of all $n\times n$
matrices over a field $\mathbb{K}$ where $n\geq 2$. If $0\leq s\leq\frac{n}{2}$, then for each $y\in M_n^{<s}({\mathbb{K}})$, there exist
$y_1,y_2\in M_n^{s}({\mathbb{K}})$ such that $y=y_1y_2$. \label{corollary-y-y1-y2}
\end{cor}

{\bf Proof.} Denote by $k\in\{0,1,\ldots,s-1\}$ the rank of $y$. From $2s\leq n$ we have $s>k\geq 0\geq 2s-n$. Then Remark \ref{remark-y-y1-y2} works.
\hfill$\Box$

The following Lemma \ref{lemma-0-any-any-0} shows  that $\delta(0)=0$ and gives  a kind  of special case for Lemma \ref{lemma-main}. Furthermore, Lemma \ref{lemma-0-any-any-0} will be used in the proof of Lemma \ref{lemma-main}.

\begin{lem}$\delta(0)=0$ for $\delta\in\mathcal{D}_s^{\times}(M_n(\mathbb{K}))$, where  $n\geq 2$ and $1\leq s\leq\frac{n}{2}$ are fixed.
In particular, for $x,y\in M_n({\mathbb{K}})$ such that $0\in \{x,y\}$, $\delta(xy)=\delta(x)y+x\delta(y)$.\label{lemma-0-any-any-0}
\end{lem}

{\bf Proof.}
Let $e=\sum_{i=1}^se_{ii}$, $f=\sum_{i=s+1}^{2s}e_{ii}$ and $g=\sum_{i=1}^se_{i,i+s}$. Certainly, $e,f,g\in  M_n^{s}({\mathbb{K}})$. From the property satisfied by $\delta$, we have
$$
\left\{
\begin{array}{l}
\delta(e)=\delta(e^2)=\delta(e)e+e\delta(e)~~\text{and} \\ \delta(f)=\delta(f^2)=\delta(f)f+f\delta(f),
\end{array}
\right.$$
which means that
$$
\left\{
\begin{array}{l}
\delta(e)=eA(I_n-e)+(I_n-e)Ae~~\text{and} \\ \delta(f)=fB(I_n-f)+(I_n-f)Bf
\end{array}
\right.
$$
for some $A,B\in M_n({\mathbb{K}})$. By the property satisfied by $\delta$, $ef=0$ and $ge=0$, we have
$$
\left\{
\begin{array}{lcl}
\delta(0)&=&\delta(ef)=e\delta(f)+\delta(e)f\\
&=&efB(I_n-f)+e(I_n-f)Bf+eA(I_n-e)f+(I_n-e)Aef=eBf+eAf~~\text{and} \\
\delta(0)&=&\delta(ge)=\delta(g)e+g\delta(e)=\delta(g)e+geA(I_n-e)+g(I_n-e)Ae=\delta(g)e+gAe,
\end{array}
\right.$$
which implies
$$\delta(0)=eBf+eAf=(eBf+eAf)f=\delta(0)f=(\delta(g)e+gAe)f=0.$$
Particularly, for any  $x\in M_n^{\leq s}({\mathbb{K}})$,
$$
\left\{
\begin{array}{l}
\delta(x\cdot0)=\delta(0)=0=x\cdot 0+\delta(x)\cdot 0=x\cdot\delta(0)+\delta(x)\cdot 0~~\text{and}\\
\delta(0\cdot x)=\delta(0)=0=0\cdot x+0\cdot\delta(x)=\delta(0)\cdot x+0\cdot\delta(x).
\end{array}
\right.
$$
\hfill$\Box$

The following lemma will be used in the proof of Theorem \ref{Theorem-rank-s}

\begin{lem}
For  $\delta\in\mathcal{D}_s^{\times}(M_n(\mathbb{K}))$ and $x,y\in M_n^{\leq s}({\mathbb{K}})$ such that $\{x,y\}\cap M_n^{\leq 1}({\mathbb{K}})\neq \phi$, where  $n\geq 2$ and $1\leq s\leq\frac{n}{2}$ are fixed, $\delta(xy)=\delta(x)y+x\delta(y)$.\label{lemma-main}
\end{lem}

{\bf Proof.}  By Lemma \ref{lemma-0-any-any-0}, it is enough to consider the case $\{x,y\}\cap M_n^{1}({\mathbb{K}})\neq \phi$. We will only prove the case for $x\in M_n^{1}({\mathbb{K}})$ and $y\in  M_n^{\leq s}({\mathbb{K}})$. The proof of the case for
$x\in M_n^{\leq s}({\mathbb{K}})$ and $y\in M_n^{1}({\mathbb{K}})$ is similar and so omitted.

{\bf Step 1.} We will prove  that for all $x\in  M_n^1({\mathbb{K}})$ and all $y\in M_n^{s}({\mathbb{K}})$, $\delta(xy)=\delta(x)y+x\delta(y)$. There exist
$P,Q,R,S\in GL_n({\mathbb{K}})$ such that $x=Pe_{11}Q$ and $y=R(\sum_{i=1}^se_{ii})S$.

{\bf Case-I.} $e_{11}QR(\sum_{i=1}^se_{ii})\neq0$. Then there exist $2\leq i_2<i_3<\cdots<i_s\leq n$ such that
 $$\left(e_{11}+\sum_{j=2}^se_{i_j,i_j}\right)QR\left(\sum_{i=1}^se_{ii}\right)\in M_n^s({\mathbb{K}})$$
since the rank of $QR(\sum_{i=1}^se_{ii})$ is $s$. By $n\geq 2s$, we can choose $2\leq k_2<k_3<\cdots<k_s\leq n$ such that
$$\{i_2,i_3,\ldots,i_s\}\cap \{k_2,k_2,\ldots,k_s\}=\phi.$$
Set $x_1=P(e_{11}+\sum_{j=2}^se_{k_j,k_j})$ and $x_2=(e_{11}+\sum_{j=2}^se_{i_j,i_j})Q$. In this case, $x_1, x_2, x_2y\in  M_n^s({\mathbb{K}})$ and $x=x_1x_2$. Hence
$$
\begin{array}{lcl}
\delta(xy)&=&\delta(x_1x_2y)=\delta(x_1(x_2y))
=\delta(x_1)x_2y+x_1\delta(x_2y)=\delta(x_1)x_2y+x_1\delta(x_2)y+x_1x_2\delta(y)\\
&=&\delta(x_1x_2)y+x_1x_2\delta(y)=\delta(x)y+x\delta(y).
\end{array}$$

{\bf Case-II.} $e_{11}QR(\sum_{i=1}^se_{ii})=0$, which means that
$$QR\left(\sum_{i=1}^se_{ii}\right)=\left(
\begin{array}{cc}
O&O\\
G&O
\end{array}\right),$$
where $G$ is an $(n-1)\times s$ matrix.
Note that the rank of $G$ is $s$. So there exist linearly independent $\alpha_1,\alpha_2,\ldots,\alpha_{n-1-s}\in \mathbb{K}^{n-1}$ such that
$$G^{T}\alpha_i=0,~~~~~~i=1,2,\ldots,n-1-s.$$
Since $2s\leq n$, we have $s-1\leq n-1-s$. Then $H=(\alpha_1,\alpha_2,\ldots,\alpha_{s-1})^{T}$ is a $(s-1)\times(n-1)$ matrix over $\mathbb{K}$. Note that the rank of $H$ is $s-1$ and $HG=0$.
Set $x_1=P(e_{11}+\sum_{i=s+1}^{2s-1}e_{ii})$ and
$$x_2=\left(
\begin{array}{cc}
1&O\\O&H\\
O&O
\end{array}
\right)Q,$$
then $x_1,x_2\in  M_n^s({\mathbb{K}})$, $x=x_1x_2$ and $x_2y=0$. Hence by the property satisfied by $\delta$ and Lemma \ref{lemma-0-any-any-0},
we have
$$
\begin{array}{lcl}
\delta(xy)&=&\delta(x_1x_2y)=\delta(x_1(x_2y))=\delta(x_1\cdot 0)
=\delta(x_1)\cdot 0+x_1\cdot\delta(0)=\delta(x_1)x_2y+x_1\delta(x_2y)\\
&=&\delta(x_1)x_2y+x_1\delta(x_2)y+x_1x_2\delta(y)=\delta(x_1x_2)y+x_1x_2\delta(y)=\delta(x)y+x\delta(y).
\end{array}
$$

{\bf Step 2.} For
 $x\in M_n^{1}({\mathbb{K}})$ and $y\in  M_n^{< s}({\mathbb{K}})$, by Corollary \ref{corollary-y-y1-y2}, there exist
$y_1,y_2\in M_n^{s}({\mathbb{K}})$ such that $y=y_1y_2$. Furthermore, by the property satisfied by $\delta$,  Lemma \ref{lemma-0-any-any-0} and the conclusion of Step 1, keeping $xy_1\in M_n^{\leq 1}({\mathbb{K}})$ in mind, we have
 $$
 \begin{array}{rcl}
 \delta(xy)&=&\delta(xy_1y_2)=\delta((xy_1)y_2)=\delta(xy_1)y_2+xy_1\delta(y_2)\\
 &=&\delta(x)y_1y_2+x\delta(y_1)y_2+xy_1\delta(y_2)=\delta(x)y_1y_2+x\delta(y_1y_2)=
 \delta(x)y+x\delta(y),
 \end{array}
 $$
which completes the proof.\hfill$\Box$

\begin{thm} Let $n$ and $s$ be integers such that $n>1$ and $1\leq s\leq \frac{n}{2}$. Let $M_n(\mathbb{K})$ be the ring of all $n\times n$ matrices over a field $\mathbb{K}$. If  a map
$\delta:M_n(\mathbb{K})\rightarrow M_n(\mathbb{K})$ satisfies that $\delta(xy)=\delta(x)y+x\delta(y)$ for any two rank-$s$ matrices $x,y$, then there exists a derivation $D$ of
$M_n(\mathbb{K})$ such that $\delta=D$ on $M_n^{\leq s}(\mathbb{K})$.\label{Theorem-rank-s}
\end{thm}

{\bf Proof.} By Lemma \ref{lemma-main}, for $i\in \underline{n}$,
\begin{equation}
\delta(e_{ii})=\delta(e_{ii}e_{ii})=e_{ii}\delta(e_{ii})+\delta(e_{ii})e_{ii},\label{equation-delta-e-ii}
\end{equation}
which means that
$\delta(e_{ii})=(a_{st}^{(i)})$, where $a_{st}^{(i)}=0$ for all
 $$(s,t)\in \left\{\left.(s,t)\in \underline{n}\times \underline{n}~\right|~s\neq i ~\text{and}~ t\neq i\}\cup \{(i,i)\right\}.$$
By Lemmas \ref{lemma-0-any-any-0} and \ref{lemma-main}, we have that for $i\neq j\in\underline{n}$
$$0=\delta(0)=\delta(e_{ii}e_{jj})=e_{ii}\delta(e_{jj})+\delta(e_{ii})e_{jj}=a_{ij}^{(j)}e_{ij}+a_{ij}^{(i)}e_{ij},$$
which means that
\begin{equation}
a_{ij}^{(j)}=-a_{ij}^{(i)} \label{equation-i-j-negative}
\end{equation}
 for $i\neq j\in\underline{n}$. Set
$$B=\sum_{i=1}^n\delta(e_{ii})e_{ii}.$$
It is easy to see that $Be_{jj}=\delta(e_{jj})e_{jj}$ for each $j\in\underline{n}$. On the other hand, by \eqref{equation-i-j-negative}
$$e_{jj}B=\sum_{k\in\underline{n}-\{j\}}a_{jk}^{(k)}e_{jk}=-\sum_{k\in\underline{n}-\{j\}}a_{jk}^{(j)}e_{jk}=-e_{jj}\delta(e_{jj}).$$
Furthermore, from Lemma \ref{lemma-main}, we have that for $i\in \underline{n}$,
$$\delta(e_{ii})=\delta(e_{ii})e_{ii}+e_{ii}\delta(e_{ii})=Be_{ii}-e_{ii}B=[B,e_{ii}].$$
For $i,j\in\underline{n}$, denote by $\lambda_{ij}$ the $(i,j)$ entry of $\delta(e_{ij})$. Note that from the expression of $\delta(e_{ii})$ we
have that $\lambda_{ii}=0$ for all $i\in\underline{n}$. For $i\neq j\in\underline{n}$, by Lemma \ref{lemma-main}
$$
\begin{array}{rcl}
\delta(e_{ij})&=&\delta(e_{ii}e_{ij}e_{jj})=\delta(e_{ii})e_{ij}+e_{ii}\delta(e_{ij})e_{jj}+e_{ij}\delta(e_{jj})\\
&=&\delta(e_{ii})e_{ii}e_{ij}+e_{ii}\delta(e_{ij})e_{jj}+e_{ij}e_{jj}\delta(e_{jj})\\
&=&Be_{ii}e_{ij}+\lambda_{ij}e_{ij}-e_{ij}e_{jj}B=[B,e_{ij}]+\lambda_{ij}e_{ij}.
\end{array}
$$
Hence for all $i,j\in\underline{n}$,
\begin{equation}
 \delta(e_{ij})=[B,e_{ij}]+\lambda_{ij}e_{ij}. \label{equation-delta-e-i-j}
\end{equation}
 For $i,j,k\in\underline{n}$, by Lemma \ref{lemma-main} and \eqref{equation-delta-e-i-j}
$$
\begin{array}{rcl}
[B,e_{ik}]+\lambda_{ik}e_{ik}&=&\delta(e_{ik})=\delta(e_{ij}e_{jk})=e_{ij}\delta(e_{jk})+\delta(e_{ij})e_{jk}\\
&=&e_{ij}([B,e_{jk}]+\lambda_{jk}e_{jk})+([B,e_{ij}]+\lambda_{ij}e_{ij})e_{jk}\\
&=&(e_{ij}[B,e_{jk}]+[B,e_{ij}]e_{jk})+(\lambda_{jk}e_{ij}e_{jk}+\lambda_{ij}e_{ij}e_{jk})\\
&=&[B,e_{ik}]+(\lambda_{ij}+\lambda_{jk})e_{ik},
\end{array}
$$
which implies that $\lambda_{ik}=\lambda_{ij}+\lambda_{jk}$ for all $i,j,k\in\underline{n}$. In particular,
$0=\lambda_{ii}=\lambda_{ij}+\lambda_{ji}$ implies that $\lambda_{ij}=-\lambda_{ji}$ for all $i, j\in\underline{n}$.
Set
$$A=B+\sum_{j=1}^n\lambda_{j1}e_{jj}.$$
Then for all $s,t\in\underline{n}$,
$$
\begin{array}{rcl}
\delta(e_{st})&=&[B,e_{st}]+\lambda_{st}e_{st}=[B,e_{st}]+(\lambda_{s1}-\lambda_{t1})e_{st}\\
&=&[B,e_{st}]+[\sum_{j=1}^n\lambda_{j1}e_{jj},e_{st}]=[A,e_{st}].
\end{array}
$$
Let $\delta'=\delta-ad_{A}$, so $\delta'\in \mathcal{D}_s^{\times}(M_n(\mathbb{K}))$ has the same property with $\delta$ and $\delta'(e_{st})=0$ for all $s,t\in\underline{n}$ by Lemma \ref{lemma-vector-space}.
For $a\in \mathbb{K}$ and $i,j\in\underline{n}$, by Lemma \ref{lemma-main},
 $$\delta'(ae_{ij})=\delta'(e_{ii}(ae_{ij})e_{jj})=\delta'(e_{ii})(ae_{ij})+e_{ii}\delta'(ae_{ij})e_{jj}+(ae_{ij})\delta'(e_{jj})=e_{ii}\delta'(ae_{ij})e_{jj},$$
 which means that there exists a map $\mu_{ij}:\mathbb{K}\rightarrow\mathbb{K}$ such that $\delta'(ae_{ij})=\mu_{ij}(a)e_{ij}$.
For $a\in\mathbb{K}$ and $i\neq j\in\underline{n}$, by Lemma \ref{lemma-main} and $\delta'(e_{st})=0$ for all $s,t\in\underline{n}$,
$$u_{ii}(a)e_{ii}=\delta'(ae_{ii})=
\left\{
\begin{array}{lcl}
\delta'((ae_{ij})e_{ji})&=&\delta'(ae_{ij})e_{ji}=u_{ij}(a)e_{ii},\\
\delta'(e_{ij}(ae_{ji}))&=&e_{ij}\delta'(ae_{ji})=u_{ji}(a)e_{ii},
\end{array}
\right.
$$
which means that $u_{ii}=u_{ij}=u_{ji}=u_{jj}$ for all $i\neq j\in\underline{n}$. Moreover for all $s,t,s',t'\in \underline{n}$,
$$u_{st}=u_{ss}=u_{s's'}=u_{s't'}.$$
Denote $u_{st}$ for $s,t\in \underline{n}$ by $\mu$.
For $a,b\in\mathbb{K}$,
$$\mu(ab)e_{11}=\delta'(abe_{11})=\delta'((ae_{11})(be_{11}))=(ae_{11})\delta'(be_{11})+\delta'(ae_{11})(be_{11})=(a\mu(b)+\mu(a)b)e_{11},$$
which implies that $\mu(ab)=a\mu(b)+\mu(a)b$ for all $a,b\in\mathbb{K}$. In particular, $\mu(1)=0$ from $\mu(1)=\mu(1\cdot 1)=\mu(1)\cdot 1+1\cdot\mu(1)$. Furthermore, for $x=(x_{ij})\in M_n^{\leq s}(\mathbb{K})$, by Lemma \ref{lemma-main} and  $\delta'(e_{st})=0$ for all $s,t\in\underline{n}$,
$$\mu(x_{ij})e_{ij}=\delta'(e_{ii}xe_{jj})=\delta'(e_{ii})xe_{jj}+e_{ii}\delta'(x)e_{jj}+e_{ii}x\delta'(e_{jj})=e_{ii}\delta'(x)e_{jj},$$
which means that $\delta'((x_{ij}))=(\mu(x_{ij}))$ for all $x=(x_{ij})\in M_n^{\leq s}(\mathbb{K})$.
Hence for all $x=(x_{ij})\in M_n^{\leq s}(\mathbb{K})$
$$\delta(x)=ad_A(x)+\delta'(x)=[A,x]+(\mu(x_{ij})).$$
At last, by Lemma \ref{lemma-main}, $\mu(1)=0$ and $\delta'((x_{ij}))=(\mu(x_{ij}))$ for all $x=(x_{ij})\in M_n^{\leq s}(\mathbb{K})$, we have that
$$
\begin{array}{rcl}
\mu(a+b)e_{11}&=&\delta'((a+b)e_{11})=\delta'((e_{11}+e_{12})(ae_{11}+be_{21}))\\
&=&(e_{11}+e_{12})\delta'(ae_{11}+be_{21})+\delta'(e_{11}+e_{12})(ae_{11}+be_{21})\\
&=&(e_{11}+e_{12})(\mu(a)e_{11}+\mu(b)e_{21})=(\mu(a)+\mu(b))e_{11},
\end{array}$$
which implies that $\mu$ is additive and so $\mu$ is a derivation of $\mathbb{K}$, further inducing a derivation $\overline{\mu}$  of $M_n(\mathbb{K})$. Note that
the restriction of  $\overline{\mu}$ on $M_n^{\leq s}(\mathbb{K})$ is $\delta'$. Hence $D=ad_{A}+\overline{\mu}$ is the desired.\hfill$\Box$

\section{Application}

For $a\in\mathbb{R}$, let $[a]$ be the least integer being not less than $a$. For example, $[-0.75]=0$, $[1.5]=2$ and $[3]=3$.
As an application of Theorem \ref{Theorem-rank-s}, we will show that the multiplicative derivation on some nonadditive subset of the matrix ring $M_n(\mathbb{K})$ over a field $\mathbb{K}$ has to be a derivation.
\begin{cor}
 Let  $n>1$ be an integer. Let $M_n(\mathbb{K})$ be the ring of all $n\times n$ matrices over a field $\mathbb{K}$. If  a map
$\delta:M_n(\mathbb{K})\rightarrow M_n(\mathbb{K})$ satisfies that $\delta(xy)=\delta(x)y+x\delta(y)$ for any
$$x,y\in\bigcup_{i=0}^{[\log_2(0.5n+1)]}M_n^{n+1-2^i}(\mathbb{K}),$$ then $\delta$ is a derivation of\label{corolary-log}
$M_n(\mathbb{K})$.
\end{cor}

{\bf Proof.} From $[\log_2(0.5n+1)]\geq \log_2(0.5n+1)$, we obtain
$$s=n+1-2^{[\log_2(0.5n+1)]}\leq n+1-2^{\log_2(0.5n+1)}=0.5n,$$
which gives that there exists a derivation $D$ of
$M_n(\mathbb{K})$ such that $\delta=D$ on $M_n^{\leq s}(\mathbb{K})$ by Theorem \ref{Theorem-rank-s}. Let $\delta'=\delta-D$, then
$\delta'(x)=0$ for all  $x\in M_n^{\leq s}(\mathbb{K})$ and $\delta'$ has the same property as $\delta$. Obviously,
 \begin{equation}
 \delta'(xy)=\delta'(x)y+x\delta'(y),~~~  x,y\in M_n^{\leq s}(\mathbb{K})\label{equation-theorem-cor-leq-s}
 \end{equation}
  and
$$
\left\{
\begin{array}{l}
\delta'(x\cdot0)=\delta'(0)=0=x\cdot 0+\delta'(x)\cdot 0=x\cdot\delta'(0)+\delta'(x)\cdot 0~~\text{and}\\
\delta'(0\cdot x)=\delta'(0)=0=0\cdot x+0\cdot\delta'(x)=\delta'(0)\cdot x+0\cdot\delta'(x),
\end{array}
\right.
$$
for all $x\in M_n(\mathbb{K})$. For any $x\in  M_n^{1}(\mathbb{K})$ and any $y\in \cup_{i=0}^{[\log_2(0.5n+1)]}M_n^{n+1-2^i}(\mathbb{K})$, there exists $x_1,x_2\in M_n^{s}(\mathbb{K})$ such that $x=x_1x_2$.
Then by $x_1,x_2,x_2y\in  M_n^{\leq s}(\mathbb{K})$ and \eqref{equation-theorem-cor-leq-s}, we have
\begin{equation}
\begin{array}{rcl}
\delta'(xy)&=&\delta'(x_1x_2y)=\delta'(x_1(x_2y))=\delta'(x_1)x_2y+x_1\delta'(x_2y)\\
&=&\delta'(x_1)x_2y+x_1\delta'(x_2)y+x_1x_2\delta'(y)=\delta'(x_1x_2)y+x_1x_2\delta'(y)\\
&=&\delta'(x)y+x\delta'(y) \label{equation-theorem-cor-1-any}
\end{array}
\end{equation}
for all $x\in M_n^{1}(\mathbb{K})$ and all $y\in \cup_{i=0}^{[\log_2(0.5n+1)]}M_n^{n+1-2^i}(\mathbb{K})$.
Similarly, we have
\begin{equation}
\delta'(yx)=\delta'(y)x+y\delta'(x) \label{equation-theorem-cor-any-1}
\end{equation}
for all $x\in M_n^{1}(\mathbb{K})$ and all $y\in \cup_{i=0}^{[\log_2(0.5n+1)]}M_n^{n+1-2^i}(\mathbb{K})$.
Hence
for all $y\in \cup_{i=0}^{[\log_2(0.5n+1)]}M_n^{n+1-2^i}(\mathbb{K})$, by \eqref{equation-theorem-cor-leq-s}, \eqref{equation-theorem-cor-1-any} and \eqref{equation-theorem-cor-any-1}, we have that for all $i,j\in\underline{n}$
$$
0=\delta'(e_{ii}ye_{jj})=\delta'(e_{ii})ye_{jj}+e_{ii}\delta'(y)e_{jj}+e_{ii}y\delta'(e_{jj})=e_{ii}\delta'(y)e_{jj}
$$
which means that $\delta'(y)=0$ for all $y\in \cup_{i=0}^{[\log_2(0.5n+1)]}M_n^{n+1-2^i}(\mathbb{K})$.

Now we only need to show that
$\delta'(z)=0$ for any rank-$k$ matrix $z$, where $n+1-2^{i+1}< k<n+1-2^i$ and $1\leq i<[\log_2(0.5n+1)]$.
It is easy to see that $$n+1-2^i>k\geq n+1-2^{i+1}+1=2(n+1-2^i)-n.$$
Then by Remark \ref{remark-y-y1-y2}, there exist $z_1,z_2\in M_n^{n+1-2^i}(\mathbb{K})$ such that $z=z_1z_2$, which implies
$$\delta'(z)=\delta'(z_1z_2)=\delta'(z_1)z_2+z_1\delta'(z_2)=0,$$
since $\delta'(y)=0$ for all $y\in \cup_{i=0}^{[\log_2(0.5n+1)]}M_n^{n+1-2^i}(\mathbb{K})$.

In conclusion, $\delta'=0$, which means that $\delta=D$ is a derivation of $M_n(\mathbb{K})$. \hfill $\Box$

% ----------------------------------------------------------------
\bibliographystyle{amsplain}
\bibliography{150}

\end{document}